\documentclass[12pt]{amsart}

\usepackage{amsfonts, amsmath}

  \usepackage{graphics}
\usepackage{amssymb}
  \usepackage{amsmath}
  \usepackage{latexsym}
\usepackage{amsthm}
\usepackage[all]{xy}

\def\C{\mathord{\mathbb C}}

\def\G{\mathord{\mathbb G}}
\def\H{\mathord{\mathbb H}}
\def\R{\mathord{\mathbb R}}
\def\K{\mathord{\mathbb K}}
\def\O{\mathord{\mathbb O}}
\def\P{\mathord{\mathbb P}}
\def\Z{\mathord{\mathbb Z}}

\def\S{\mathord{\mathbb S}}

\def\V{{\mathcal V}}
\def\M{{\mathcal M}}
\def\A{{\mathcal A}}

\def\3{{\ss}}
\def\2{\frac{1}{2}}
\def\nhalf{\begin{smallmatrix} \underline{n} \cr 2 \end{smallmatrix}}
\def\x{\times}
\def\.{\cdot}
\def\d{\partial}
\def\o{\circ}

\def\<{\langle}
\def\>{\rangle}

\def\Fix{\mathop{\rm Fix}\nolimits}
\def\diag{\mathop{\rm diag}\nolimits}

\def\Span{\mathop{\rm Span\,}\nolimits}

\def\trace{\mathop{\rm trace\,}\nolimits}

\def\so{\mathfrak{so}}
\def\id{{\rm id}}
\def\i{{\sf i}}

\def\ms{\medskip\noindent}

\def\beq{\begin{equation}}
\def\eeq{\end{equation}}
\def\bea{\begin{eqnarray}}
\def\eea{\end{eqnarray}}
\def\bsm{\left(\begin{smallmatrix}}
\def\esm{\end{smallmatrix}\right)}
\def\bpm{\begin{pmatrix}}
\def\epm{\end{pmatrix}}

\begin{document}
  \title{Bott Periodicity, Submanifolds, and Vector Bundles}
  
  \author{Jost Eschenburg}
  \author{Bernhard Hanke}
  \address{Institut f\"ur Mathematik, Universit\"at Augsburg, D-86135 Augsburg, Germany}
  \email{eschenburg@math.uni-augsburg.de}
  \email{hanke@math.uni-augsburg.de}
  

 \date{\today}


\maketitle

\newtheorem{theorem}{\bf Theorem}
\newtheorem{corollary}[theorem]{\bf Corollary}
\newtheorem{proposition}[theorem]{\bf Proposition}
\newtheorem{lemma}[theorem]{\bf Lemma}
\newtheorem{definition}[theorem]{\bf Definition}

\begin{abstract}{We sketch a geometric proof of the classical theorem of Atiyah, Bott, and Shapiro \cite{ABS}
    which relates Clifford modules to vector bundles over spheres. 
    Every module of the Clifford algebra $Cl_k$ defines a particular vector bundle over $\S^{k+1}$,
    a generalized Hopf bundle, and the theorem asserts that this correspondence
    between $Cl_k$-modules and stable vector bundles over $\S^{k+1}$ is an isomorphism
    modulo $Cl_{k+1}$-modules. We prove this theorem directly, 
    based on explicit deformations as in Milnor's book on Morse theory \cite{M}, and
    without referring to the Bott periodicity theorem as in \cite{ABS}.
}
\end{abstract}

\section*{Introduction}

Topology and Geometry are related in various ways. Often topological properties
of a specific space are obtained by assembling its local curvature invariants,
like in the Gauss-Bonnet theorem. Bott's periodicity theorem is different: A detailed
investigation of certain totally geodesic submanifolds in specific symmetric spaces
leads to fundamental insight not just for these spaces but for whole areas
of mathematics. This geometric approach was used originally by Bott \cite{Bott1,Bott2} 
and Milnor in his book on Morse theory \cite{M} where the stable homotopy 
of the classical groups was computed. Later Bott's periodicity theorem 
was re-interpreted as a theorem on K-theory \cite{AB,ABS,A}, 
but the proofs were different and less geometric. 
However we feel that the original approach of Bott and Milnor 
can prove also the K-theoretic versions of the periodicity theorem. 
As an example we discuss Theorem (11.5.) from the fundamental paper \cite{ABS} 
by Atiyah, Bott and Shapiro, which relates
Clifford modules to vector bundles over spheres. 
The argument in \cite{ABS} uses explicit computations
of the right and left hand sides of the stated isomorphism, 
and depends on the Bott periodicity theorem
for the orthogonal groups. Instead we prove bijectivity of the 
relevant comparison map directly. In
consequence the Bott periodicity theorem for the orthogonal groups
is now implied by its algebraic counterpart in the representation 
theory of Clifford algebras \cite{ABS}. This gives a positive response to the
remark in \cite[page 4]{ABS}: ``It is to be hoped that Theorem (11.5) 
can be give a more natural and less computational
proof'', cf.\ also \cite[page 69]{LM}. We will concentrate
on the real case which is more interesting and less well known than the
complex theory. Much of the necessary geometry was explained to us
by Peter Quast \cite{Q1}.

\section{Poles and Centrioles}

We start with the geometry. A {\it symmetric space} is a Riemannian manifold $P$
with an isometric point reflection $s_p$ (called {\it symmetry}) at any point $p\in P$, 
that is $s_p \in \hat G$ = isometry group of $P$ with $s_p(\exp_p(v)) = \exp_p(-v)$ 
for all $v\in T_pP$. The map $s : p\mapsto s_p : P \to \hat G$ is called 
{\it Cartan map}; it is a covering onto its 
image $s(P)\subset\hat G$ which is also symmetric.\footnote 
	{$s(P) \subset \hat G$ is a connected component of the set $\{g\in\hat G: g^{-1} = g\}$.
	When we choose a symmetric metric on $\hat G$ such that $g\mapsto g^{-1}$ is an isometry,
	$s(P)$ is a reflective submanifold and hence totally geodesic, thus symmetric.}
The composition of any
two symmetries, $\tau = s_q s_p$ is called a {\it transvection}. It translates the 
geodesic $\gamma$ connecting $p = \gamma(0)$ to $q = \gamma(r)$ by $2r$ and acts by parallel 
translation along $\gamma$, see next figure. 
The subgroup of $\hat G$ generated by all transvections (acting transitively on $P$) will be called $G$.

\ms\hskip .7cm
\includegraphics{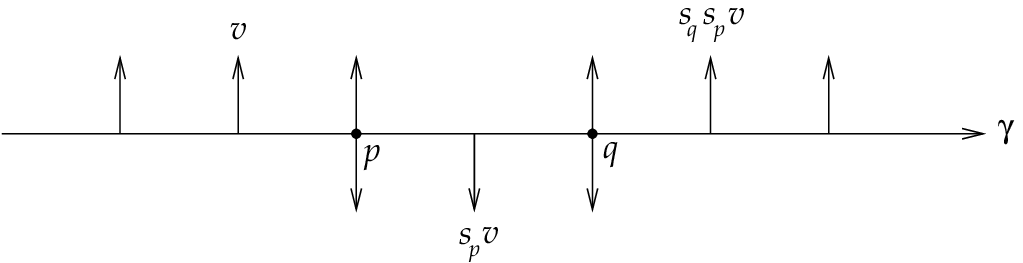}

\ms Two points $o,p\in P$ will be called {\it poles}  if $s_p = s_o$. 
The notion was coined for the north and south pole of a round sphere, but there
are many other spaces with poles; e.g.\ $P = SO_{2n}$ with $o = I$ and 
$p = -I$, or the Grassmannian $P=\G_n(\R^{2n})$ with $o = \R^n$ and $p = (\R^n)^\perp$. 
A geodesic $\gamma$ connecting $o = \gamma(0)$ to $p = \gamma(1)$ 
is reflected into itself at $o$ and $p$ and hence it is closed with period $2$.

\ms Now we consider the {\it midpoint set} $M$ between poles $o$ and $p$,
$$
  M = \{m=\gamma\bsm \underline 1 \cr 2 \esm: \gamma \textrm{ geodesic in $P$ with } 
	\gamma(0) = o,\ \gamma(1) = p\}.
$$
For the sphere
$P=\S^n$ with north pole $o$, this set would be the equator, see figure below. 

\hskip 3.5cm
\includegraphics{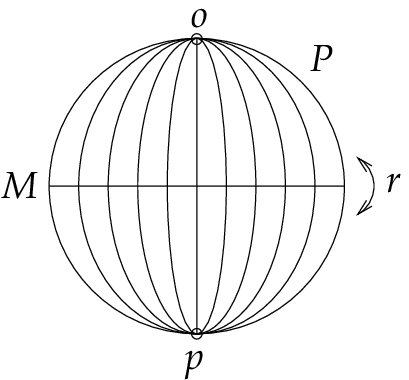}

\begin{theorem} {\rm \cite{N}} $M$ is the fixed set of an isometric involution $r$ on $P$.
\end{theorem}
 
\proof In the example of the sphere $P = \S^n$, the equator $M$ is the fixed set of $-s_o = -I\o s_o$.
Here, $-I$ is the deck transformation\footnote 
	{A {\it deck transformation} of  
	$\pi : P\to \bar P$ is an isometry $\delta$ of $P$ with $\pi\o\delta = \pi$.}
of the covering $\S^n \to \R\P^n = \S^n/\{\pm I\}$. In the general case
we consider the covering $P \to s(P)$. Since $s(P)$ is again symmetric, we have
$s(P) = P/\Delta$ for some discrete freely acting group $\Delta\subset \hat G$
normalized by all symmetries and centralized by all transvections.\footnote
	{Consider a symmetric space $P$ and a covering $\pi:P\to P/\Delta$ 
	for some discrete freely acting group $\Delta$ of isometries
	on $P$. Then $P/\Delta$ is again symmetric if and only if each symmetry $s_p$ of $P$ maps $\Delta$-orbits
	onto $\Delta$-orbits. Thus for each $\delta\in\Delta$ we have $s_p(\delta x) = \tilde\delta s_p(x)$ 
	for all $x\in P$, and $\tilde\delta\in\Delta$
	is independent of $x$, by discreteness. Thus $s_p\delta = \tilde\delta s_p$, in particular
	$s_p\delta s_p = \tilde\delta\in\Delta$. For any other symmetry $s_q$ we have the same equation
	$s_q\delta = \tilde\delta s_q$ with the same $\tilde\delta\in\Delta$, again by discreteness. Thus
	$\delta^{-1}s_ps_q\delta = s_p\tilde\delta^{-1}\tilde\delta s_q = s_ps_q$, and $\delta$ commutes with the
	transvection $s_ps_q$ (see also \cite[Thm.\ 8.3.11]{W}).}
Since $s_o = s_p$, the points $o$ and $p$ are identified in $s(P)$.  
Thus there is a unique $\delta \in \Delta$ with $\delta(o) = p$. 
This will be the analogue of $-I$ in the case $P = \S^n$.
We will show that $\delta$ has order 2 and preserves any geodesic $\gamma$ 
with $\gamma(0) = o$ and $\gamma(1) = p$. In fact, let $\tau$ 
be the transvection along $\gamma$ from $o$ to $p$. Then $\tau^2(o) = o$ and therefore
$$
	\delta(p) = \delta(\tau(o)) = \tau(\delta(o)) = \tau(p) = o.
$$
Thus $\delta^2$ fixes $o$ which shows $\delta^2 = \id$ since $\Delta$ acts freely.
Hence $\{I,\delta\}\subset\Delta$ is a subgroup and $\bar P = P/\{\id,\delta\}$ a
symmetric space. Under the projection $\pi : P \to \bar P$, the geodesic $\gamma$ is
mapped onto a closed geodesic doubly covered by $\gamma$, thus $\delta$ preserves $\gamma$
and shifts its parameter by $1$, and $\gamma$ has period 2.

\ms
\hskip 4.3cm
\includegraphics{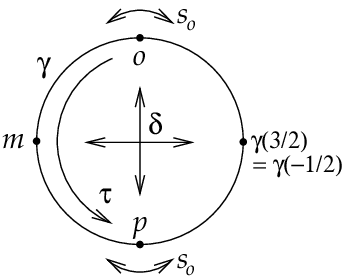}

\ms We put $r = s_o\delta$. This is an involution since $s_o$ and $\delta$ commute:
$\delta'= s_o\delta s_o\in\Delta$ sends $o$ to $p$ like $\delta$, thus $\delta'=\delta$. 
Then $r$ fixes the midpoint $m=\gamma(\frac{1}{2})$ of any geodesic
$\gamma$ from $o$ to $p$ since $s_o(\delta(\gamma(\frac{1}{2})) = s_o(\gamma(\frac{3}{2}))
= s_o(\gamma(-\2)) =  \gamma(\frac{1}{2})$. Thus $M \subset \Fix(r)$.

Vice versa, assume that $m\in P$ is a fixed point of $r$. Thus $s_om = \delta m$. Join $o$ to $m$ by a geodesic 
$\gamma$ with $\gamma(0) = o$ and $\gamma(\2) = m$. Then $\gamma(-\2) = s_o(m) = \delta(m) = \delta(\gamma(\2))$, 
and the projection $\pi : P \to \bar P = P/\{\id,\delta\}$ maps $\gamma : [-\2,\2] \to P$ 
onto a geodesic loop $\bar\gamma = \pi\o\gamma$, that is a closed geodesic of period 1 (since $\bar P$ is symmetric). 
Thus $\gamma$ extends to a closed geodesic of period 2 doubly covering $\bar\gamma$, and $\delta$ shifts the
parameter of $\gamma$ by 1. Therefore $\gamma(1) = \delta(o) = p$. Hence $m$ is the midpoint of $\gamma|_{[0,1]}$
from $o$ to $p$. Thus $M \supset \Fix(r)$.
\endproof

Connected components of the midpoint set $M$ are called {\it centrioles} \cite{CN}.
Connected components of the fixed set of an isometry are totally geodesic 
(otherwise shortest geodesic segments in the ambient space with end points in the fixed set 
were not unique, see figure below); if the isometry is an involution, its fixed
components are called {\it reflective}.

\ms
\hskip 3.5cm
\includegraphics{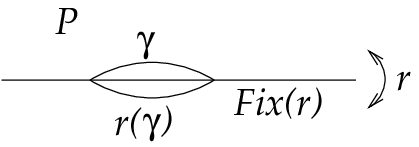}

\ms Most interesting are connected components containing midpoints of geod\-esics with {\it minimal} length
between $o$ and $p$ (``{\it minimal centrioles}''). 
Each such midpoint $m = \gamma(\2)$ determines its geodesic $\gamma$
uniquely: if there were two geodesics of equal length from $o$ to $p$ through $m$, they
could be made shorter by cutting the corner. \label{corner}

\ms\hskip 3.5cm
\includegraphics{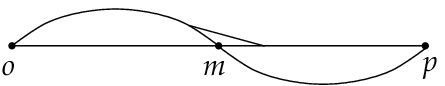}

\medskip There exist chains of minimal centrioles (centrioles in centrioles),
\beq \label{P}
	P \supset P_1 \supset P_2 \supset \dots 
\eeq 
Peter Quast \cite{Q1,Q2} classified all such chains with at least 3 steps
starting with a compact simple Lie group $P = G$. Up to group coverings, the result is as follows.
The chains 1,2,3 occur in Milnor \cite{M}.
\beq \label{table}
\begin{matrix}
\textrm{No.} & G & P_1 & P_2 & P_3 & P_4 & \textrm{restr.}	\cr
\hline
1 & (S)O_{4n} & SO_{4n}/U_{4n} & U_{2n}/Sp_n & \G_p(\H^n) & Sp_p & p=\nhalf \cr
2 & (S)U_{2n} & \G_n(\C^{2n}) & U_n & \G_p(\C^n) & U_p & p=\nhalf \cr
3 & Sp_n & Sp_n/U_n & U_n/SO_n & \G_p(\R^n) & SO_p & p=\nhalf \cr
4 & Spin_n & Q_n & (\S^1\!\x\!\S^{n-3})/\!\pm &\S^{n-4} & \S^{n-5} & n\geq5 \cr
5 & E_7 & E_7/(\S^1E_6) & \S^1E_6/F_4 & \O\P^2 & -
\end{matrix}
\eeq

\ms 
By $\G_p(\K^n)$ we denote the Grassmannian of $p$-dimensional subspaces in
$\K^n$ for $\K \in \{\R,\C,\H\}$. Further, $Q_n$ denotes the complex quadric in $\C\P^{n+1}$ 
which is isomorphic to the real Grassmannian $\G^+_2(\R^{n+2})$ of oriented 2-planes, 
and $\O\P^2$ is the octonionic projective plane $F_4/Spin_9$.

A chain is extendible beyond $P_k$ if and only if
$P_k$ contains poles again. E.g.\ among the Grassmannians $P_3=\G_p(\K^n)$ 
only those of half dimensional subspaces ($p = \frac{n}2$)
enjoy this property: Then $(E,E^\perp)$ is a pair of poles for any $E\in\G_{n/2}(\K^n)$, 
and the corresponding midpoint set is the
group $O_{n/2}, U_{n/2}, Sp_{n/2}$ since its elements are the graphs of orthogonal 
$\K$-linear maps $E\to E^\perp$, see figure below.

\hskip 3.5cm
\includegraphics{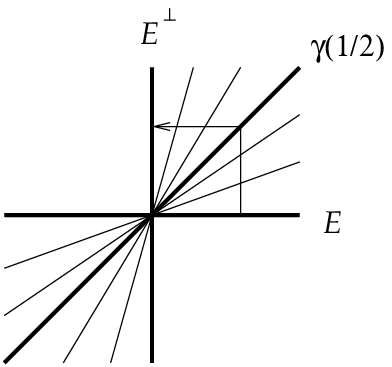}

\section{Centrioles with topological meaning} \label{topology}

Points in minimal centrioles are in 1:1 correspondence to minimal geodesics between the
corresponding poles $o$ and $p$. 
Thus minimal centrioles sometimes can be viewed as low-dimensional approximations 
of the full {\it path space} $\Lambda$,
the space of all $H^1$-curves\footnote 
	{$H^1$ means that $\lambda$ has a derivative almost everywhere which is square integrable.
	Replacing any path $\lambda$ by a geodesic polygon with $N$ vertices, we may replace
	$\Lambda$ by a finite dimensional manifold, cf.\ \cite{M}.}
$\lambda : [0,1] \to P$ with $\lambda(0) = o$ and $\lambda(1) = p$. 
This is due to the Morse theory for the energy function $E$ on $\Lambda$ where
$E(\lambda) = \int_0^1 |\lambda'(t)|^2dt$. We may decrease the energy of
any path $\lambda$ by applying the gradient flow of $-E$ (left figure).

\ms\hskip 1.5cm
\includegraphics{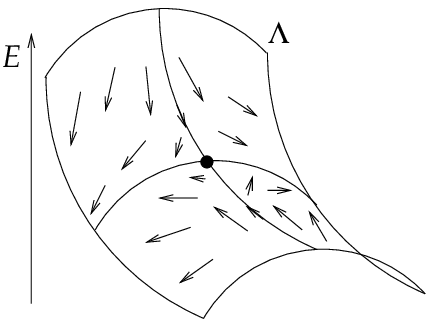} \hskip 1cm
\includegraphics{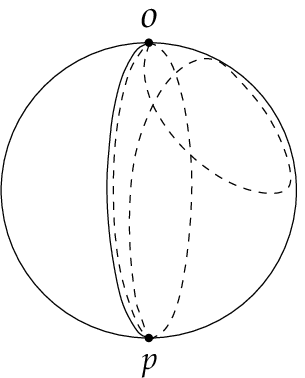}

\ms
Most elements of $\Lambda$
will be flowed to the minima of $E$ which are the shortest geodesics between $o$ and $p$. 
The only exceptions are the domains of
attraction (``unstable manifolds'') for the other critical points, the non-minimal geodesics between $o$ and $p$. 
The codimension of the unstable manifold is the {\it index} of the critical point, the maximal dimension
of any subspace where the second derivative of $E$ (taken at the critical point) is negative. 
If $\beta$ denotes the smallest index of all non-minimal critical points, any continuous map 
$f : X \to \Lambda$ from a connected cell complex $X$ of dimension $<\beta$ can be moved away from these
unstable manifolds and flowed into a connected component of the minimum set, that is into some centriole $P_1$.
Thus $f$ is homotopic to a map $\tilde f : X\to P_1$.

But this works only if all non-minimal geodesics from $o$ to $p$ have high index ($\geq\beta$). 
Which symmetric spaces $P$ have this property? An easy example is the sphere, $P = \S^n$.
A nonminimal geodesic $\gamma$ between poles $o$ and $p$ covers a great circle at least one and a half times
and can be shortened within any 2-sphere in which it lies (right figure above). There are $n-1$ such 2-spheres
perpendicular to each other since the tangent vector $\gamma'(0)
= e_1$ is contained in $n-1$ perpendicular planes $\Span(e_1,e_i)$ with $i\geq2$ in the tangent space. 
Thus the index is $\geq n-1$, in fact $\geq 2(n-1)$ since any such geodesic contains at least 2 conjugate
points where it can be shortened by cutting the corner, see figure.

\ms\hskip2.3cm
\includegraphics{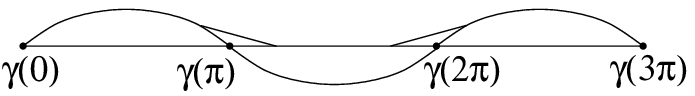}

\ms
For the classical groups we can argue similarly. E.g.\ in $SO_{2n}$, a shortest geodesic from $I$ to $-I$
is a product of $n$ half turns, planar rotations by the angle $\pi$ in $n$ perpendicular 2-planes in $\R^{2n}$. 
A non-minimal geodesic must make an additional full turn and thus a $3\pi$-rotation in at least one 
of these planes, say in the $x_1x_2$-plane. This rotation belongs to the rotation group 
$SO_3 \subset SO_{2n}$ in the $x_1x_2x_k$-space for any $k\in\{3,\dots,2n\}$. Using $SO_3 = \S^3/\pm$,
we lift the $3\pi$-rotation to $\S^3$ and obtain a 3/4 great circle which can be shortened. 
There are $2n-2$ coordinates $x_k$ and therefore $2n-2$ independent contracting directions, hence the
index of a nonminimal geodesic in $SO_{2n}$ is $\geq2n-2$ (compare \cite[Lemma 24.2]{M}).
The index of the spaces $P_k$ can be bounded from below in a similar way, see next section for
the chain of $SO_n$.
This implies the homotopy version of the periodicity theorem:

\begin{theorem} \label{bott}
When $n$ is even and sufficiently large, we have for $G = SO_{4n},SU_{2n},Sp_n$ (notations of table \ref{table}):
$$
	\pi_k(G) = \pi_{k-1}(P_1)= \pi_{k-2}(P_2) = \pi_{k-3}(P_3)) = \pi_{k-4}(P_4) .
$$
Together with table \ref{table} this implies the following periodicities: 
\bea
	\pi_{k+2}(SU_n) &=& \pi_k(SU_{n/2}),\cr 
	\pi_{k+4}(SO_{n}) &=& \pi_k(Sp_{n/8}),\cr
	\pi_{k+4}(Sp_n) &=& \pi_k(SO_{n/2}).\nonumber
\eea
\end{theorem}

\section{Clifford modules} \label{clifford}

For compact matrix groups $G$ containing $-I$, there is a linear algebra interpretation for the iterated 
midpoint sets $M_j$ and their components $P_j$. A geodesic $\gamma$ in $G$ with $\gamma(0) = I$ is a one-parameter
subgroup, and when $\gamma(1) = -I$, then $\gamma(\2) = J$ is a {\it complex structure}, $J^2 = -I$. 
Thus the midpoint set $M_1$ is the set of complex structures in $G$. When the connected component $P_1$
of $M_1$ contains antipodal points $J_1$ and $-J_1$, there is a next midpoint set $M_2\subset P_1$. It 
consists of points $J_1\gamma(\2)$ where $\gamma$ is a one-parameter subgroup in $G$ with $\gamma(1)=-I$ such that $J_1\gamma(t)$ is a complex structure for all $t$, 
$$
	J_1\gamma J_1\gamma=-I\,.
\eqno{(*)}
$$ 
In particular the midpoint $J = \gamma(\2)$ 
anticommutes with $J_1$ (since $J_1JJ_1J = -I$ $\iff$ $J_1J = -JJ_1$), and when $\gamma$ is minimal, this
condition is sufficient for ($*$): 
then both $J_1\gamma J_1$ and $-\gamma^{-1}$ are shortest geodesics from $-I$ to $I$
with midpoint $J$, so they must agree. By induction hypothesis, we have anticommuting complex structures
$J_u\in G$ with $J_i\in P_i$ for $i<k$, and $P_{k}$ is a connected component of the set 
\beq \label{Pk}
	M_k = \{J\in G: J^2 = -I,\, JJ_i = -J_iJ \textrm{ for } i<k\}
\eeq
of complex structures $J\in G$ which anticommute with $J_1,\dots,J_{k-1}$. 
To finish the induction step we choose some $J_k \in P_k$. 

Recall that the {\it real Clifford algebra} $Cl_k$ is the associative real algebra with 1 which is generated 
by $\R^k$ with the relations $vw+wv = -2\<v,w\>$. Equivalently, an orthonormal basis $e_1,\dots,e_k$ of $\R^k
\subset Cl_k$ satisfies 
$$
	e_ie_j+e_je_i = -2\delta_{ij}\,.
$$ 
A {\it representation}
of $Cl_k$ is an algebra homomorphism from $Cl_k$ into some matrix algebra $\K^{n\x n}$ with $\K\in\{\R,\C,\H\}$;
the space $\K^n$ on which the matrices operate is called {\it Clifford module} $S$. A representation maps
the vectors $e_i$ onto matrices $J_i$ with the same relations $J_i^2 = -I$ and $J_iJ_j = -J_jJ_i$ for $i\neq j$.
Thus a $Cl_k$ module is nothing but a {\it Clifford system}, a family of $k$ are anticommuting complex structures, 
and the midpoint set $M_{k+1}\subset P_k$ between $J_k$ and $-J_k$ can be viewed as the set of 
extensions of a given $Cl_{k}$-module (defined by $J_1,\dots,J_k$) to a $Cl_{k+1}$-module. 

\medskip The algebraic theory of the Clifford representations is rather easy (cf.\ \cite{LM}). 
They are direct sums of
irreducible representations, and in the real case there is just one 
irreducible $Cl_k$-module $S_k$ (up to isomorphisms) when $k \not\equiv 3\mod 4$, while
there are two with equal dimensions when $k \equiv 3 \mod 4$.  
For $k = 0,\dots,8$ we have 

\begin{theorem}
\beq \label{Mk}
\begin{matrix}
	k & 0 & 1 & 2 & 3 & 4 & 5 & 6 & 7 & 8 	\cr
	S_k & \R &\C&\H&\H&\H^2&\C^4&\O&\O&\O^2
\end{matrix}
\eeq 
and further we have the {\it periodicity theorem for Clifford modules},
\beq \label{Mk+8}
	S_{k+8} = S_k \otimes S_8.
\eeq
For $k=3$ and $k=7$, the two different module structures are given by 
left and right multiplications of $\R^k=\K':= \K\ominus \R\.1$ on $S_k=\K$ for $\K = \H,\O$. 
\end{theorem}

\section{Index of nonminimal geodesics} \label{index}

From (\ref{Pk}) we have gained a uniform description for all iterated centrioles $P_k$ of $G$ 
in terms of Clifford systems.
This can be used for a calculation of the lower bound $\beta$ 
for the index of nonminimal geodesics in all $P_k$, cf.\ \cite{M}.\footnote
	{A different argument using root systems was given by Bott[(6.7)] \cite{Bott1} and in more
	detail by Mitchell \cite{Mt1,Mt2}}

\begin{theorem} Let $SO_n = G \supset P_1 \supset P_2 \supset\dots\supset P_k\supset\dots$ be the chain (\ref{P}) 
of itererated centrioles where $n$ is divisible by a high power of $2$. Then for each $k$ there is some lower bound 
$\beta$ depending on $n$ such that the index of nonminimal geodesics from $J_k$ to $-J_k$ is $\geq\beta$, and $\beta \to \infty$ as $n\to \infty$.
\end{theorem}

\proof
Let $\tilde\gamma = J_k\gamma: [0,1]\to P_k$ be a non-minimal geodesic from $J_k$ to $-J_k$. Then 
$\gamma(t) = e^{\pi tA}$ for some $A\in\so_{n}$. 
Since $\tilde\gamma(t)$ anticommutes with $J_i$ for all $i<k$, it follows that 
$\gamma(t)$ and $A$ commute with $J_i$. Further, from $\tilde\gamma(t)^2 = -I$ we obtain 
$J_{k} e^{\pi tA} J_{k}^{-1} = e^{-\pi tA}$ and therefore $A$ anticommutes with $J_{k}$.
Thus we have computed the tangent space of $P_k$ at $J_k$:
\beq \label{TPk}
	T_{J_k}P_k = \{J_kA: A\in \so_n,\ AJ_k = -J_kA,\ AJ_i = J_iA \textrm{ for } i<k\}.
\eeq
Since $\gamma(1) = -I$, the (complex) eigenvalues of $A$ 
have the form $a\i$ with $\i = \sqrt{-1}$ and $a$ an odd integer.

To relate these eigenvalues to the index we argue similar as in \cite[p.\ 144-147]{M}. 
We split $\R^{n}$ into a sum of subspaces $V_j$ being invariant
under the linear maps $A,J_1,\dots,J_{k}$ and being minimal with respect to this property. 
All $J_i$, $i<k$, preserve the (complex) eigenspaces $E_a$ of $A$,
corresponding to the nonzero eigenvalue $a\i$, while $J_{k}$ interchanges $E_a$ and $E_{-a}$.
Thus by minimality of $V_j$, there is just one pair $\pm a$ such that $V_j$ is the real part
of $E_a+E_{-a}$. Therefore $J' := A/a$ is an additional complex structure on $V_j$ commuting with 
$J_i$ ($i<k$) 
and anticommuting with $J_{k}$, and $J_{k+1} := J_{k}J'$ is a complex structure which 
anticommutes with {\em all} $J_1,\dots,J_{k}$. Hence 
$V_j$ is an irreducible $Cl_{k}$-module.
Moreover, $A = a_jJ'$ on $V_j$ for some nonzero
integer $a_j$ while $A = 0$ on $V_0$. By choice of the sign of $J'|V_j$ we may assume that all $a_j > 0$. 
hence $a_j \in\{1,3,5,\dots\}$.

Choose two of these irreducible modules, say $V_j$ and $V_h$. 
By (\ref{Mk}), there is a module isomorphism $V_j \to V_h$
as $Cl_{k+1}$-modules when $k+1 \not\equiv 3 \mod 4$ (Case 1) and as $Cl_k$-modules when $k+1\equiv 3\mod 4$
(Case 2). 
This remains true when we alter the $Cl_{k+1}$-module structure of $V_h$ 
in Case 1 by changing the sign of $J_{k+1}$ (and thus that of $J'$) on $V_h$. 
With this identification we have $V_j+V_h = V_j\otimes\R^2$ 
and $B = I\otimes\bsm &-I\cr I \esm$ (with $B=0$ on the other submodules) 
commutes with all $J_j$, $j\leq k$, and the same is true for $e^{uB}$. 
Putting $A_u = e^{uB}A e^{-uB}$, we have $J_{k}A_u\in T_{J_{k}}P_k$ by (\ref{TPk}).

\ms {\bf Case 1:} $\underline{k+1 \not\equiv 3 \mod 4}$: 
We have modified our identification of $V_j$ and $V_h$ by changing the sign of $J_{k+1}$ on $V_h$.
Thus on $V_j+V_h = V_j\otimes\R^2$ we have $A = J' \otimes D$  where 
	$D = \diag(a_j,-a_h) = cI+D'$
with $D' = \diag(b,-b)$ for $b = \2(a_j+a_h)$ and $c = \2(a_j-a_h)$.  
Let us consider the family of geodesics $J_k\gamma_u$ from $J_k$ to $-J_k$ in $P_k$ with 
$\gamma_u(t) = e^{t\pi A_u} = e^{uB}\gamma(t)e^{-uB}$.
The point $\gamma(t) = e^{\pi tc}e^{\pi tD'}$
is fixed under conjugation with the rotation matrix $e^{uB} = \bsm \cos u & -\sin u \cr \sin u & \cos u \esm$
precisely when $e^{\pi tD'} = \diag( e^{\pi tb} , e^{-\pi tb})$ is a multiple of the identity matrix
which happens for $t =1/b$. If one of the eigenvalues $a$ of $A$
is $>1$, say $a_h \geq 3$, then $b =\2(a_j+a_h) \geq 2$ and $1/b \in (0,1)$.
All $\gamma_u$ are geodesics connecting $I$ to $-I$ on $[0,1]$. 
By ``cutting the corner'' it follows that $\gamma$ can no longer be locally
shortest beyond $t\!=\! 1/b$ , see figure. If there is at least one eigenvalue $a_h > 1$,
we have $r-1$ index pairs $(j,h)$, hence the index of non-minimal geodesics is at least $r-1$.

\ms
\hskip 3.5cm
\includegraphics{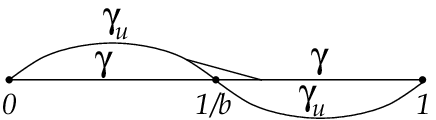}

\ms {\bf Case 2:} $\underline{k+1 \equiv 3 \mod 4}$:
In this case, the product $J_o := J_1J_2\dots J_{k-1}$ is a complex structure\footnote
	{Putting $S_n = (J_1\dots J_n)^2$ we have 
	$$S_n = J_1\dots J_nJ_1\dots J_n = (-1)^{n-1}S_{n-1}J_n^2
	= (-1)^nS_{n-1},$$
	thus $S_n = (-1)^sI$ with $s = n+(n-1)+\dots+1 = \2 n(n+1)$. When $n = k-1 \equiv 1\mod 4$,
	then $s$ is odd, hence $S_n = -I$.}
which commutes with $A$ and anticommutes with $J_k$ (since $k-1$ is odd). 
Thus $A$ can be viewed as a complex matrix,
using $J_o$ as the multiplicaton by $\i$. 
Let $E_a\subset V_j$ be the eigenspace of $A$ corresponding to the eigenvalue $\i a$ where $a$ is any
odd integer. Then $E_a$ is
invariant under the $J_i$, $i<k$, which commute with $A$, but is it also invariant under $J_k$ which
anticommutes with $A$ and with $\i = J_o$ (since $k-1$ is odd). By minimality we have $V_j = E_a$, hence
$A = aJ_o$. As before, we consider the linear map $B = \bsm &-I\cr I\esm$ on $V_j+V_h = V_j\otimes \R^2$
and the family of geodesics $\gamma_u(t) = e^{t\pi A_u} = e^{uB}\gamma(t)e^{-uB}$. This time, 
$A = J' \otimes D$  where $D = \diag(a_j,a_h) = cI+D'$
with $c = \2(a_j+a_h)$ and $D' = \diag(b,-b)$ with $b = \2(a_j-a_h)$. Thus the element 
$\gamma(t) = e^{\pi tc}e^{\pi tD'}$
is fixed under conjugation with the rotation matrix $e^{uB} = \bsm \cos u & -\sin u \cr \sin u & \cos u \esm$
precisely when $e^{\pi tD'} = \diag( e^{\pi tb} , e^{-\pi tb})$ is a multiple of the identity matrix
which happens for $t =1/b$. If $b>1$, we obtain an energy-decreasing deformation by cutting the corner, 
see figure above. We need to show that there are enough pairs $(j,h)$ with $b > 1$ when $\gamma$ is non-minimal.

Any $J \in P_k$ defines a $\C$-linear map $J_kJ$ since $J_kJ$ commutes with $J_i$ and hence with $J_o$.
Thus a path $\lambda : I \to P_k$ from $J_k$ to $-J_k$ defines a family of $\C$-linear maps, and its
complex determinant $\det(J_k\lambda)$ is a path in $\S^1$ starting and ending at $\det(\pm I) = 1$
(recall that the dimension $n$ is even).
This loop in $\S^1$ has a mapping degree which is apparently invariant under homotopy; it decomposes the path
space $\Lambda P_k$ into infinitely many connected components. 
If $\lambda$ is a geodesic, $\lambda(t) = J_k e^{\pi t A}$,
then $\det J_k^{-1}\lambda(t) = e^{\pi t\,\trace A}$, hence its mapping degree is $\2\trace A/\i$.
Since $\trace A/\i = m\sum_j a_j$, we may fix $c := \sum_j a_j$ 
(which means fixing the connected component of $\Lambda P_k$) and we may
assume that $|c|$ is much smaller than $r$ (the number of submodules $V_j$). Let $p$ denote the
sum of the positive $a_j$ and $-q$ the sum of the negative $a_j$. Then $p+q \geq r$ since all $|a_j|\geq 0$,
and $p-q = c$ which means roughly $p\approx q \approx r/2$. Assume for the moment $c = 0$. 
If there is some eigenvalue $a_h$ with $|a_h| > 1$, say $a_h = -3$, there are many positive $a_j$
with $a_j-a_h\geq 4$, more precisely $\sum_{a_j>0} (a_j-a_h) \geq 4\.r/2 = 2r$,
and this is a lower bound for the index. In the general case this result has to be corrected by the
comparably small number $c$. In contrast, if all $a_j = \pm 1$, the geodesic $\gamma$ consists of simultaneous
half turns in $n/2$ perpendicular planes; these are shortest geodesics from $I$ to $-I$ in $SO_n$.\footnote
	{Any one-parameter subgroup $\gamma$ in $SO_n$ is a family of planar rotations in perpendicular planes.
	When $\gamma(1) = -I$, all rotation angles are odd multiples of $\pi$. The squared length of $\gamma$
	is the sum of the squared rotation angles. Thus the length is minimal if all rotation angles
	are just $\pm \pi$.}

\section{Vector bundles over spheres}

Clifford representations have a direct connection to vector bundles over spheres
and hence to K-theory. Every vector bundle $E \to \S^{k+1}$ is trivial 
over each of the two closed hemispheres 
$D_+,D_-\subset \S^{k+1}$, but along the equator $\S^k = D_+\cap D_-$ the fibres over $\d D_+$ and $\d D_-$ are identified by some map $\phi:\S^k \to O_n$ called {\it clutching map}. 

\ms\hskip 3.5cm
\includegraphics{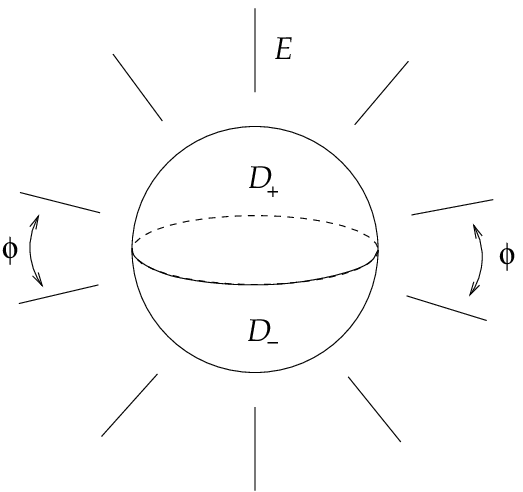}

\ms Homotopic clutching maps define equivalent
vector bundles. Thus vector bundles over $\S^{k+1}$ are classified by the homotopy group $\pi_k(O_n)$. When we
allow adding of trivial bundles (stabilization), $n$ may be arbitrarily high. Let $\V_k$ be the set of
vector bundles over $\S^{k+1}$ up to equivalence and adding of trivial bundles (``stable vector bundles''). Then
\beq \label{Vk=pik}
	\V_k = \lim_{n\to\infty} \pi_k(O_n).
\eeq
Hence we could apply Theorem \ref{bott} in order to classify stable vector bundles over spheres. 
However, a seperate argument based on the same ideas but also using Clifford modules will give more information.

\medskip A $Cl_k$ module $S=\R^n$ or the corresponding Clifford system $J_1,\dots,J_k$ $\in O_n$ defines a 
peculiar map $\phi = \phi_S: \S^k\to O_n$ which is {\it linear}, that is a restriction of a linear map 
$\phi:\R^{k+1}\to\R^{n\x n}$, where we put
\beq
	\phi_S(e_{k+1}) = I,\ \ \ \phi_S(e_i) = J_i \textrm{ for } i\leq k .
\eeq
The bundles defined by such clutching maps $\phi_S$ are called {\it generalized Hopf bundles}. In the cases
$k=1,3,7$, these are the classical complex, quaternionic, and octonionic Hopf bundles over $\S^{k+1}$.

In fact, $Cl_k$-modules are in 1:1 correspondence to linear maps $\phi:\S^k\to O_n$ with the identity
matrix in the image. To see this, let $\phi$ be such map and $W = \phi(\R^{k+1})$ its image. Then $\S_W := 
\phi(S^k)\subset O_n$. Thus $\phi$ is an isometry for the inner product 
$\<A,B\> = \frac{1}{n}\trace(A^TB)$ on $\R^{n\x n}$ since $\phi(\S^k) \subset O_n$ 
and $O_n$ lies in the unit sphere of $\R^{n\x n}$. 
For all $A,B\in \S_W$ we have $(A+B) \in \R\.O_n$. On the other hand, 
$(A+B)^T(A+B) = 2I + A^TB+B^TA$, thus $A^TB + B^TA = tI$ 
for some $t\in\R$. From the inner product with $I$ we obtain $t = 2\<A,B\>$. Inserting $A = I$ and $B\perp I$ 
yields $B+B^T = 0$, and for any $A,B\perp I$ we obtain $AB+BA = -2\<A,B\>I$. Thus 
$\phi|\R^k$ defines a $Cl_k$-representation on $\R^n$.

\medskip
Atiyah, Bott and Shipiro \cite{ABS} reduced the theory of vector bundles
over spheres to the simple algebraic structure of Clifford modules by showing that all vector bundles
over spheres are generalized Hopf bundles plus trivial bundles, see Theorem \ref{ABS} below. 
We sketch a different proof of this theorem using the original ideas of Bott and Milnor. 
We will homotopically deform the clutching map $\phi : \S^k \to G = SO_{n}$ of
the given bundle $E \to \S^{k+1}$ step by step into a linear map. Since adding of trivial bundles
is allowed, we may assume that the rank $n$ of $E$ is divisible by a high power of 2.

\medskip We declare $N = e_{k+1}$ to be the ``north pole'' of $\S^k$. First we deform $\phi$ such
that $\phi(N) = I$ and $\phi(-N) = -I$. Thus $\phi$ maps each meridian from $N$ to $-N$ in $\S^k$ 
onto some path from $I$ to $-I$ in $G$, an element of $\Lambda G$. 
The meridians $\mu_v$ are parametrized by $v\in \S^{k-1}$ 
where $\S^{k-1}$ is the equator of $\S^k$. Therefore $\phi$ can be considered as a map $\phi : \S^{k-1}\to\Lambda G$.
Using the negative gradient flow for the energy function $E$ on the path space $\Lambda G$ 
as in section \ref{topology} 
we may shorten all $\phi(\mu_v)$ simultaneously to minimal geodesics from $I$ to $-I$ and obtain a map
$\tilde\phi : \S^{k-1} \to \Lambda_oG$ where $\Lambda_oG$ is the set of shortest
geodesics from $I$ to $-I$, the minimum set of $E$ on $\Lambda G$. 
Let $m(\gamma) = \gamma(\2)$ be the midpoint of any geodesic $\gamma :
[0,1]\to G$. Thus we obtain a map $\phi_1 = m\o\tilde\phi : \S^{k-1}\to P_1$, and we may replace $\phi$
by the geodesic suspension over $\phi_1$ from $I$ and $-I$.

\ms\hskip .5cm
\includegraphics{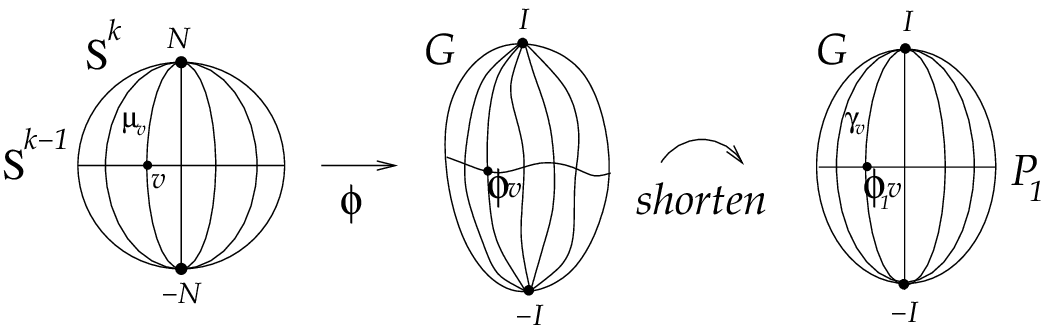} 

\ms
We repeat this step replacing $G$ by $P_1$ and $\phi$ by $\phi_1$.
Again we choose a ``north pole'' $N_1 = e_k \in \S^{k-1}$ and deform $\phi_1$ such that $\phi_1(\pm N_1)
= \pm J_1$ for some $J_1\in P_1$. Now we deform the curves $\phi_1(\mu_1)$ for all meridians $\mu_1\subset\S^{k-1}$ to  shortest geodesics, whose midpoints define a map $\phi_2 : \S^{k-2} \to P_1$, and then we replace $\phi_1$ by a geodesic suspension from $\pm J_1$ over $\phi_2$. 
This step is repeated $(k-1)$-times until we reach a map $\phi_{k-1} :
\S^1\to P_{k-1}$. This loop can be shortened to a geodesic 
loop $\tilde\gamma = J_{k-1}\gamma : [0,1] \to P_{k-1}$ (which is a closed geodesic since $P_{k-1}$ is symmetric) 
starting and ending at $J_{k-1}$, such that $\tilde\gamma$ and $\gamma$ are shortest in their homotopy class.

We have $\gamma(t) = e^{2\pi tA}$ for some $A\in T_{J_{k-1}}P_{k-1}$.  
Since $\gamma$ is closed, the (complex) eigenvalues of $A$ 
have the form $a\i$ with $a\in\Z$ and $\i = \sqrt{-1}$.
To compute these eigenvalues we argue as in section \ref{index}. 
We split $\R^{n}$ into $V_0 = \ker A$ and a sum of subspaces $V_j$ which are invariant
under the linear maps $A,J_1,\dots,J_{k-1}$ and minimal with respect to this property. 
As before, $A = aJ'$ for some nonnegative integer $a$, and $J_k = J_{k-1}J'$ is a complex structure 
anticommuting with $J_1,\dots,J_{k-1}$. Hence $V_j$ is an irreducible $Cl_{k}$-module with dimension $m_k$,
see (\ref{Mk}), (\ref{Mk+8}). 

Since the clutching map of the
given vector bundle $E\to \S^{k+1}$ (after the deformation) is determined by $\gamma,J_1,\dots,J_{k-1}$ 
which leave all $V_j$, $j\geq 0$, invariant, the vector bundle splits accordingly as $E = E_0 \oplus \sum_{j>0} E_j$ 
where $E_0$ is trivial.\footnote
	{The clutching map $\phi : \S^k\to SO_{n}$ splits into components $\phi_j : \S^k \to SO(V_j)$. 
	The domain $\S^k$ is the union of totally geodesic spherical $(k-1)$-discs 
	$D_v$, $v\in\S^1$, centered at $v$ and perpendicular to $\S^1$. 
	All $D_v$ have a common boundary $\S^{k-2}$. 
	Since $\phi_0|_{D_v}$ is constant in $v$, it is contractable along $D_v$
	to a constant map.

\hskip 4cm
\includegraphics{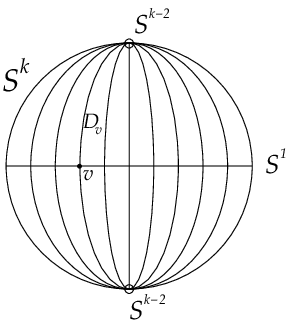}
}

 We claim that the minimality of $\gamma$ implies $a_j = 1$ for all $j$ and hence $A = J_k$.
In fact, the geodesic variaton $\gamma_u$ of section \ref{index} shows that $|a_j-a_h| < 2$ for all $j,h$,
otherwise we could shorten $\gamma$ by cutting the corner.

\ms
\hskip 3.5cm
\includegraphics{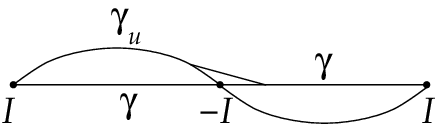}

\ms
Now suppose that, say, $a_1 \geq 2$. We may assume that $V_0=\ker A$ contains another copy $\tilde V_1$ of $V_1$
as a $Cl_{k-1}$-module: if not, we extend $E_0$ by the trivial bundle $\S^{k+1}\x \tilde V_1$.
Thus we have eigenvalues $0$ and $a_1$ on $\tilde V_1\oplus V_1$ with difference $\geq 2$, in contradiction
to the minimality of the geodesic. 

We have shown $E = E_0\oplus E_1$ 
where $E_0$ is trivial and $E_1$ is a generalized
Hopf bundle for the Clifford system $J_1,\dots,J_k$ on $\sum_{j>0} V_j$.

\medskip
Let $\M_k$ the set of equivalence classes of $Cl_k$-modules, modulo trivial $Cl_k$-representations. 
We have studied the map
$$
	\hat\alpha : \M_k \to \V_k
$$
which assigns to each $S \in \M_k$ the corresponding generalized Hopf bundle over $\S^{k+1}$.
It is additive with respect to direct sums. We have just proved that $\hat\alpha$ is onto. But it not 1:1. 
In fact, every $Cl_{k+1}$-module is also a $Cl_k$ module since $Cl_k \subset Cl_{k+1}$.
This defines a restriction map $\rho : \M_{k+1} \to \M_k$. Any $Cl_k$-module $S$ which is really
a $Cl_{k+1}$-module gives rise to a contractible clutching map $\phi_S : \S^k \to SO_{n}$ and hence
to a trivial vector bundle since $\phi_S$ can be 
extended to $\S^{k+1}$ and thus contracted over one of the half spheres $D_+,D_-\subset \S^{k+1}$. 
Thus $\hat\alpha$ sends $\rho(\M_{k+1})$ into trivial bundles and hence it descends to an additive map\footnote
	{In fact, both $\V_k$ and $\A_k$ are abelian groups with respect to direct sums, not just 
	semigroups, and $\alpha$ is a group homomorphism. 
	Using the tensor product, $\V = \sum_k \V_k$ and $\A = \sum_k \A_k$ become rings and 
	$\alpha$ a ring homomorphism, see \cite{ABS}.}
$$
	\alpha : \A_k := \M_k/\rho(\M_{k+1}) \to \V_k.
$$
We claim that this map is injective: if a bundle $\hat\alpha(S)$ is
trivial for some $Cl_k$-module $S$, then $S$ is (the restriction of) a $Cl_{k+1}$-module.

\ms
{\it Proof of the claim.} Let $S$ be a $Cl_k$-module and 
$\phi = \phi_S : \S^k \to G$ the corresponding clutching map (that is $\phi(e_{k+1}) = I$, $\phi(e_i) = J_i$).
We assume that $\phi$ is contractible, that is it extends to $\hat\phi :D^{k+1}\to G$. The closed disk
$D^{k+1}$ will be considered as the northern hemisphere $D^{k+1}_+ \subset \S^{k+1}$. 
Repeating the argument above for the surjectivity, we consider the meridians $\mu_v$ between $N = e_{k+1}\in\S^k$
and $-N$, but this time there are much more such meridians, not only those in $\S^k$ but also those through
the hemisphere $D^{k+1}_+$. They are labeled by $v \in D^k_+ := D^{k+1}_+\cap N^\perp$. 

\ms\hskip 4.5cm
\includegraphics{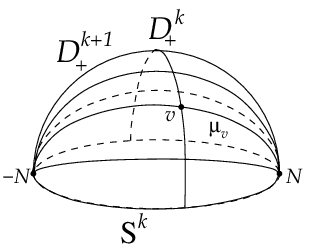}

\ms Applying the negative energy gradient flow we deform the curves $\phi(\mu_v)$ to minimal geodesics 
without changing those in $\phi(\S^k)$ which are already minimal. Then we obtain the midpoint map 
$\hat\phi_1 : D^k_+ \to P_1$ with $\phi_1(v) = m(\hat\phi(\mu_v))$ which extends the given midpoint map 
$\phi_1$ of $\phi$. This step is repeated $k$ times until we reach $\hat\phi_k : D^1_+ \to P_k$
which is a path from $J_k$ to $-J_k$ in $P_k$. This path can be shortened to
a minimal geodesic in $P_k$ whose midpoint is a complex structure $J_{k+1}$ anticommuting with 
$J_1,\dots,J_k$. Thus we have shown that our $Cl_k$-module $S$ is extendible to a $Cl_{k+1}$-module,
that is $S \in \rho(\M_{k+1})$. This finishes the proof of the injectivity.

\begin{theorem} \label{ABS} {\rm \cite{ABS}} Every vector bundle over $\S^k$ splits into a trivial bundle and a
generalized Hopf bundle. More precisely, the map
$
	\alpha : \A_k = \M_k/\rho(\M_{k+1}) \to \V_k
$
sending the equivalence class of a $Cl_k$-module $S$ onto its generalized Hopf bundle is an isomorphism.
\end{theorem}

From (\ref{Mk}) one easily obtains the groups $\A_k$ since the modules $S_k$ in (\ref{Mk}) 
are the (one or two) generators of $\M_k$. If $S_k = \rho(S_{k+1})$, then $\A_k = 0$. This happens for $k= 2,4,5,6$.
For $k = 0,1$ we have 
	$$\rho(S_{k+1}) = S_k\oplus S_k = 2S_k,$$ 
hence $\A_0 = \A_1 = \Z_2$.
For $k = 3,7$ there are two generators for $\M_k$, say $S_k$ and $S'_k$, and $\rho(S_{k+1})
= S_k \oplus S_k'$, thus $\A_3 = \A_7 = \Z$. Hence 
\beq \label{Ak}
\begin{matrix}
	k &| & 0 & 1 & 2 & 3 & 4 & 5 & 6 & 7  	\cr
	\A_k & | & \Z_2 & \Z_2 & 0 & \Z & 0 & 0 & 0 & \Z
\end{matrix}
\eeq
and because of the periodicity (\ref{Mk+8}) we have $\A_{k+8} = \A_k$.

\smallskip Consequently, the list (\ref{Ak}) for $\A_k$ is the same as that for $\V_k$ and for $\pi_k(O_n)$, $n$ large (see (\ref{Vk=pik}). Thus we have also computed the stable homotopy of $O_n$.

\medskip
We have seen that the following objects are closely related and obey
the same periodicity theorem:
\begin{itemize}
\item Iterated centrioles of $O_n$,
\item stable homotopy groups of $O_n$,
\item Clifford modules,
\item stable vector bundles over spheres.
\end{itemize}


\begin{thebibliography} {9999}

\bibitem{A} M.\ Atiyah: K-Theory. Benjamin 1967
\bibitem{AB} M.\ Atiyah, R.\ Bott, {\em On the periodicity theorem for complex vector bundles},
Acta Math.\ 112 (1964), 229 - 247
\bibitem{ABS} M.\ Atiyah, R.\ Bott, A.\ Shapiro, {\em Clifford modules}, Topology 3 (1964), 3-38
\bibitem{Bott1} R.\ Bott, {\em The stable homotopy of the classical groups}, Proc. Nat. Acad. Sci. U.S.A. {\bf 43} (1957), 933--935. 
\bibitem{Bott2} R.\ Bott, {\em The stable homotopy of the classical groups}, Ann. Math. {\bf 70} (1959), 313--337. 
\bibitem{CN} Chen, B.-Y., Nagano, T., {\em Totally geodesic submanifolds of symmetric
spaces II}, Duke Math.\ J.\ {\bf 45} (1978), 405 - 425

\bibitem{LM} H.B.\ Lawson, M.-L.\ Michelson: Spin Geometry, Princeton 1989

\bibitem{M} J.\ Milnor: Morse Theory, Princeton 1963
\bibitem{Mt1} Mitchell, S.\ A., {\it The Bott filtration of a loop group}, in: Springer Lect.\ Notes Math.\
 1298 (1987), 215 - 226

\bibitem{Mt2} Mitchell, S.\ A.: {\it Quillen’s theorem on buildings and the loops on a symmetric space}, 
Enseign.\ Math. (II)34 (1988), 123 - 166

\bibitem{N} Nagano, T., {\it The involutions of compact symmetric spaces}, Tokyo J.\ Math.
11 (1988), 57 - 79


\bibitem{Q1} P.\ Quast: {\it Complex structures and chains of symmetric
spaces}, Habilitations\-schrift Augsburg 2010 (available from the author)

\bibitem{Q2} P.\ Quast: {\it Centrioles in Symmetric Spaces}, Nagoya Math.\ J.\ 211 (2013), 51 - 77


\bibitem{W} Wolf, J.A.: {\em Spaces of constant curvature}, 5th edition, Publish or Perish,
Wilmington (Delaware) 1984

\end{thebibliography}
\end{document}